\documentclass[11pt]{amsart}

\usepackage{amsfonts,amsmath,amssymb,amsthm}
\usepackage{fancyhdr}
\usepackage[all]{xypic}
\usepackage{array}

\newcommand{\Cc}{\mathbb{C}}  
\newcommand{\Rr}{\mathbb{R}}

\renewcommand{\epsilon}{\varepsilon}

\renewcommand{\ge}{\geqslant}

\renewcommand {\geq}{\geqslant}

\newcommand{\m}{\mathfrak{m}_0}

\newcommand{\grad}{\mathop{\mathrm{grad}}\nolimits}
\newcommand{\Int}{\mathop{\mathrm{Int}}\nolimits}
\newcommand{\Baff}{{\mathcal{B}_\mathit{\!aff}}}
\newcommand{\Binf}{{\mathcal{B}_\infty}}
\newcommand{\B}{{\mathcal{B}}}
\newcommand{\Fgen}{\mathcal{F}_\mathit{\!gen}}
\newcommand{\Fgenb}{\bar{\mathcal{F}}_{\mathit{\!gen}}}
\newcommand{\zero}{ \left\{ 0 \right\} }
\newcommand{\F}{\mathcal{F}}
\newcommand{\C}{\mathcal{C}}
\newcommand{\OO}{\mathcal{O}}

\newcommand{\red}{{\mathit{\!red}}}
    
\newcommand{\Ker}{\mathop{\mathrm{Ker}}\nolimits}
\newcommand{\id}{\mathop{\mathrm{id}}\nolimits}

\renewcommand{\Im}{\mathop{\mathrm{Im}}\nolimits}

\newtheorem*{theorem*}{Theorem}  
\newtheorem{theorem}{Theorem}    
\newtheorem{lemma}[theorem]{Lemma}        
\newtheorem{proposition}[theorem]{Proposition}
\newtheorem*{remark*}{Remark}  

\begin{document}
\title{\textbf{Classification of polynomials from 
$\Cc^2$ to $\Cc$ with one critical value}}
\author{Arnaud Bodin}

\maketitle

\section{Introduction}

Let $f : \Cc^2 \longrightarrow \Cc$ be a polynomial map. The 
\emph{bifurcation set} $\B$ is the minimal set of points of $\Cc$
such that $f : \Cc^2 \setminus f^{-1}(\B) \longrightarrow \Cc \setminus \B$ is
a locally trivial fibration. 
We can describe $\B$ as follows: let 
$\Baff = \left\{ f(x,y)\ | \ \grad_f(x,y)=(0,0) \right\}$ 
be the set of \emph{affine critical values}.
The set $\Baff$ is a subset of $\B$ but is not necessarily equal to $\B$. 
The value $c\in \Cc$ is \emph{regular at infinity} if there
exists a disk ${D}$ centered at $c$ and a compact set ${K}$
of $\Cc^2$ with a locally trivial fibration $f : f^{-1} \left( {D} 
\right)\setminus {K} \longrightarrow {D}$. There is only a finite number of non-regular
values at infinity: the \emph{critical values at infinity} collected
in $\Binf$. The bifurcation set $\B$ is now:
$$\B = \Baff \cup \Binf.$$
For $c\in \Cc$, we denote the fiber $f^{-1}(c)$ by $\F_c$.
If $s\notin \B$, then the fiber $\F_s$ is  called a \emph{generic fiber} and is 
denoted $\Fgen$. 

The aim of this paper is to describe the classification of reduced poly\-no\-mial maps
with one critical value, that is, for convenience, $\B = \zero$.
The  classification is  given up to homeomorphisms: two polynomials
$f$ and $g$ are \emph{topologically equivalent} ($f \approx g$) if there exists 
homeomorphisms $\Phi$ and $\Psi$ such that the following 
diagram commutes:
$$
\xymatrix{
{}\Cc^2\ar[r]^-\Phi\ar[d]_-f  & \Cc^2\ar[d]^-g \\
\Cc \ar[r]_-\Psi        & \Cc.
}
$$

\begin{theorem*}
Let $f : \Cc^2 \longrightarrow \Cc$  be a reduced polynomial.
We denote by $p$ and $q$ two relatively prime natural numbers, 
 $\epsilon,\epsilon' \in \{0,1\}$, $\sigma = \sigma(x,y) = x^sy+1$, 
($s > 0$).
Let $n\geqslant 1$ and let $g(x)$ be the polynomial $g(x)= \prod_{i=1}^{n} (x-i)^{m_i}$ with 
$1 \leqslant m_1\leqslant m_2\leqslant \cdots \leqslant  m_n$;
 and let $g_\red$ be the reduced polynomial associated to $g$, $g_\red(x)= \prod_{i=1}^{n} (x-i)$.
\begin{itemize}
  \item If $\Baff = \Binf = \varnothing$ then $f \approx x$,
  \item if $\Baff = \zero$ and $\Binf = \varnothing$ then
  \begin{itemize}
      \item $f \approx y \cdot g_\red(x)$,
      \item or  $f \approx x\prod_{i=1}^{n}({x^p-i y})$ (if $p=1$ then 
$n \geqslant 2$),
      \item  or $f \approx {x^\epsilon y^{\epsilon'}\prod_{i=1}^{n}({x^p-i y^q})}$, ($1<p<q$),
  \end{itemize}
  \item if $\Baff = \varnothing$ and $\Binf = \zero$ then 
  \begin{itemize}
        \item $f \approx x \prod_{i=1}^{n}\left({x^py^q-i}\right)$ 
        \item or $f \approx x \sigma \prod_{i=1}^{n}\left({x^p\sigma^q-i}\right)$ 
              ($p>1$ or $q>1$),
        \item  or  $f \approx x\sigma^\epsilon \prod_{i=1}^{n}\left({ x^p-i\sigma^q}\right)$  
(if $\epsilon = 0$ then $q > 1$),
        \item or $f \approx g_\red(x)( g(x)y+1 )$ ($n > 1$),
  \end{itemize}
  \item if $\Baff = \zero$ and $\Binf = \zero$ then 
  \begin{itemize}
     \item $f \approx xy \prod_{i=1}^{n}(x^py^q-i)$ ($1\leqslant p<q$),
     \item or $f \approx g_\red(x)k(x)( g(x)y+1 )$ 
          ($k(x) = \prod_{i=1}^{n'} (x+i)$, $n'\geqslant 1$).
   \end{itemize}
\end{itemize}
Moreover, two different polynomials of this classification are not topo\-lo\-gically
equivalent.
\end{theorem*}

We  define a stronger notion of equivalence.
Two polynomials $f$ and $g$ are \emph{algebraically equivalent} ($f \sim g$) if there exists 
an algebraic automorphism $\Phi$ of $\Cc^2$ and $\Psi$ an 
automorphism  of $\Cc$ with equation $\Psi(z) = az+b$
such that the following  diagram commutes:
$$
\xymatrix{
{}\Cc^2\ar[r]^-\Phi\ar[d]_-f  & \Cc^2\ar[d]^-g \\
\Cc \ar[r]_-\Psi        & \Cc.
}
$$

The study of polynomials with one critical value is reduced
to a few cases, up to algebraic equivalence.
The case $\B=\varnothing$ is the famous Abhyankar-Moh theorem
(\cite{AM}, see paragraph \ref{sec:pre}). 
A theorem of  M.~Za\u{\i}denberg and V.~Lin \cite{ZL} 
corresponds to the case $\Baff = \zero$ and $\Binf = \varnothing$
for irreducible polynomials.
We generalize this result to the reducible case by using methods 
from the proof of Za\u{\i}denberg-Lin theorem by 
W.~Neumann and L.~Rudolph \cite{NR} (paragraph \ref{sec:ZL}).
The remaining cases ($\Binf = \zero$)  are studied in 
paragraphs \ref{sec:first} and \ref{sec:second}. 
The arguments are essentially topological: we find a smooth disk 
in the fiber $f^{-1}(0)$ and we argue with branched coverings 
in order to give equations that represent equivalent classes of polynomials up
to algebraic equivalence. That enables us to recover the list obtained 
by M.~Za\u{\i}denberg  by the use of $\Cc^*$-action \cite{Z}.

The last part of the work (paragraph \ref{sec:topo}) is to deduce from
the former results the topological classification.
Resolution of singularities determines  polynomials
with one critical value up to topological equivalence.
It gives  a classification without redundancy.
The algebraic and the topological classification for
irreducible polynomials with $\Binf = \varnothing$ (and with
$\Baff = \varnothing$ or $\Baff = \zero$) given by Abhyankar-Moh
and Za\u{\i}denberg-Lin are the same. However this is not true in general:
we give polynomials (with $\Baff = \varnothing$ and $\Binf = \zero$) that are
topologically equivalent but not algebraically equivalent.

\section{Preliminaries}
\label{sec:pre}

When there is no critical value, the situation has been completed by 
S.~Abhyankar and T.~Moh \cite{AM}.
Abhyankar-Moh theorem   is formulated as follows:
\begin{theorem}
\label{th:AM}
If $\B=\varnothing$ then $f \sim x$.
\end{theorem}

Recall that $\F_0 = f^{-1}(0)$.
A polynomial is \emph{primitive} if its generic fiber is connected.
The link between Euler characteristic of the zero fiber
and the inclusion $\B \subset \zero$ (that is to say $\B=\varnothing$ or
$\B= \zero$) is explained in the lemma:
\begin{lemma}
\label{lem:chi}
If $\B \subset \zero$ then $\chi\left( \F_0 \right) = +1$.
Moreover, if the polynomial $f$ is primitive and  
$\chi\left( \F_0 \right) = +1$ then $\B \subset \zero$.
\end{lemma}
\begin{proof}
The decomposition $\Cc = \Cc \setminus \{ 0 \} \cup \{ 0 \}$
gives a partition $\Cc^2 = f^{-1}( \Cc \setminus \{ 0 \}) \cup f^{-1}(0)$.
By additivity of the Euler characteristic, \cite[p.~95]{Fu}
$$ 1 = \chi \big(  f^{-1}( \Cc \setminus \{ 0 \}) \big) + 
\chi \big(  f^{-1}(0) \big). $$
If $\B \subset \{ 0\}$ then $f$ defines a locally trivial 
fibration onto $ \Cc \setminus \{ 0 \}$. Then
$$ \chi \big(  f^{-1}( \Cc \setminus \{ 0 \}) \big) = 
\chi \big(  \Cc \setminus \{ 0 \} \big)
\times \chi \big(  f^{-1}(1) \big).$$
The  Euler characteristic of $ \Cc \setminus \{ 0 \}$ 
is zero. Hence $ \chi \big(  f^{-1}( \Cc \setminus \{ 0 \}) \big) = 0$ and
$\chi \big( f^{-1}(0) \big) = 1 $.

Conversely, if $f$ is a primitive polynomial then by Suzuki 
 formula \cite{S}:
$$1-\chi(\Fgen) = \sum_{c\in \B}\big(\chi(\F_c)
  -\chi\big(\Fgen)\big).$$
If $\chi(\F_0)=+1$ then 
 $\sum_{c\in \B \setminus \{ 0 \} }\big(\chi(\F_c) -\chi(\Fgen)\big)=0$,
but if $c \in \B$ then $\chi(\F_c) -\chi(\Fgen) > 0$ (see \cite{HL}), then 
$\B \subset \{ 0 \}$.
\end{proof}

\begin{remark*}
For a primitive polynomial Suzuki formula proves the equivalence $\chi \big( f^{-1}(0) \big) = +1 
\Leftrightarrow \B \subset \zero$ (see \cite{ZL}, \cite{GP} for example). However the non 
primitive polynomial $f(x,y) = xy(xy+1)$ verifies $\chi \big( f^{-1}(0) \big) = 1$ but 
$\B = \{ 0, -\frac{1}{4} \}$.
\end{remark*}

\bigskip

We denote $h(0)$ the algebraic monodromy induced in 
homology\footnote{Homology with integer coefficients.} on
$H_1(\Fgen)$
by a small circle $S^1_\epsilon(0)$ of radius $\epsilon$ centered  at $0$.
The key of this paper is the following simple remark: for all $S^1_r(0)$ ($r>0$)
the induced monodromies are equal since $0$ is the only critical value.

To compactify the situation we need resolution of singularities at infinity
\cite{LW}:
$$
\xymatrix{
{}\Cc^2   \ar[r] \ar[d]_-f      &\Cc P^2 \ar[d]_-{\tilde{f}}    
    &\Sigma \ar[l]_-{\pi} \ar[dl]^-{\bar{f}}   \\
\Cc \ar[r]      &\Cc P^1
}
$$
where $\tilde{f}$ is the natural ---but not well-defined--- map
coming from the homogenization of $f$; $\pi$ is the blow-up
of some points on the line at infinity $L_\infty$ of $\Cc P^2$ and 
of the affine singular points.

We denote $D_0 = \bar{f}^{-1}(0)$ and $D_\infty= \bar{f}^{-1}(\infty)$.
The \emph{dual graph} $G_0$ of $D_0$ is obtained as follows:
one vertex for each irreducible component of $D_0$ and one edge between two
vertices for one intersection of the corresponding components. 
A similar construction is done to obtain $G_\infty$, we know
 that $G_\infty$ is a tree \cite{LW}.

The monodromy induced by a small circle $S^1_\epsilon(\infty)$  centered
 at $\infty$ in $\Cc P^1$ is exactly the monodromy $h(0)$ with the 
reverse orientation:
$$ h(\infty) = h(0)^{-1}.$$
This property allows us to prove the three following lemmas.
\begin{lemma}
\label{lem:rationnal}
The fiber $\F_0=f^{-1}(0)$ is rational, that is to say the union 
of punctured spheres.
\end{lemma}
\begin{proof}
Let $B_1,\ldots,B_p$ be small $4$-balls around the affine
singularities of $\F_0$ and set 
$\F_0^\circ = \F_0 \cap B^4_R \setminus B_1\cup\ldots\cup B_p$.
Then $\F_0^\circ$ can be isotoped into $\Fgen$ and we denote
$\ell_*:H_1(\F_0^\circ) \longrightarrow H_1(\Fgen)$ the induced morphism.
Then, by \cite{B} or \cite{MW}, the invariant cycles for $h(0)$ are
$\Ker(h(0)-\id)= \Im \ell_*$.
Suppose that one of the components of $\F_0$ has genus,
then $\Fgen$ has genus and the
cycles corresponding to genus induced be $\F_0^\circ$ are invariant cycles.

On the other hand the cycles invariant by all the monodromies
associated to elements of $\pi_1(\Cc\setminus \B,*)$
are cycles corresponding to the 
boundary\footnote{If $\Fgenb$ is the surface without boundary associated to $\Fgen$,
$\iota : \Fgen \longrightarrow \Fgenb$ is the inclusion and
$\iota_* : H_1(\Fgen) \longrightarrow H_1(\Fgenb)$ is the induced
morphism, then the ``boundary cycles'' are $\Ker \iota_*$.} 
of $\Fgen$ (see \cite{B} or \cite{DN}).
Here there is only one monodromy and invariant cycles
by all the monodromies are exactly cycles invariant
by $h(0)$. It provides a contradiction.
\end{proof}

\begin{lemma}
\label{lem:nocycle}
There is no cycle in $G_0$: $H_1(G_0) = 0$.
\end{lemma}
\begin{proof}
A theorem of F.~Michel and C.~Weber \cite{MW} asserts firstly, that the cycles of $G_0$ 
correspond to Jordan $2$-blocks 
$\left(\begin{smallmatrix} 1&1\\0&1 \end{smallmatrix}\right)$
 for the monodromy $h(0)$ and secondly, that
$h(\infty)$ does not have any such blocks since $G_\infty$ is a tree.
Now, as $h(0)=h(\infty)^{-1}$, $G_0$ has no cycle.
\end{proof}

\begin{lemma}
\label{lem:tubseif}
The tube $f^{-1}(S^1_r(0))$ is a Seifert manifold.
\end{lemma}

\begin{proof}
Let us suppose that in the minimal Waldhausen 
decomposition\footnote{Or the Jaco-Shalen-Johannson decomposition.}
of $f^{-1}(S^1_r(0))$  there exists two  distinct Seifert pieces.
This decomposition can be obtained, as described in \cite{LMW},
from the boundary of a neighborhood of the divisor $D_0$; moreover 
a Dehn twist between two Seifert pieces can be calculated
(see \cite{MW}) and is non-positive. But the decomposition of $f^{-1}(S^1_r(0))$
can also be obtained as the boundary of a neighborhood of $D_\infty$
(because $\B \subset \zero$). Then the same formula proves that the Dehn twist is
non-negative since the orientation of the second boundary is the opposite
of the first;  now the Dehn twist is non-negative and non-positive, hence equal to zero.
That contradicts the fact that the essential pieces were distinct.
\end{proof}

In other words, let us call a singularity that provides a Seifert piece
in the decomposition of $f^{-1}(S^1_r(0))$ an \emph{essential singularity}.
We have proved that at most one essential singularity can occur.
The non-essential affine singularities are  ordinary quadratic singularities.

\section{Generalization of Za\u{\i}denberg-Lin theorem}
\label{sec:ZL}

Let us recall Za\u{\i}denberg-Lin theorem \cite{ZL}. For a proof
using topological arguments see \cite{NR}.
\begin{theorem}
\label{th:ZL}
Let $f$ be an irreducible polynomial with the fiber $\F_0=f^{-1}(0)$ simply connected
then for some relatively prime natural numbers $p$ and $q$ (or for $p=1$ and $q=0$):
$$ f\sim x^p-y^q.$$
\end{theorem}

The following lemma links the topology of $\F_0$ to the case  without
critical value at infinity.

\begin{lemma}
\label{lem:simp}
 Let $f$ be a reduced polynomial. Then $\Baff \subset \zero$
and $\Binf=\varnothing$ if and only if $\F_0$ is a simply connected set.
\end{lemma}
\begin{proof}
If $\F_0$ is simply connected then the irregular fiber is connected and it is
reduced because $f$ is a reduced polynomial,
hence the generic fiber is also connected, then $f$ is a primitive polynomial.
Moreover $\chi(\F_0)=+1$ so $H_1(\F_0)=\zero$ and by lemma \ref{lem:chi}, 
$\B \subset \zero$. Let $T_0$ be the tube $f^{-1}(\Delta)$ where $\Delta$ is
a small\footnote{Small enough in comparison to $R$ that defines 
the link at infinity $f^{-1}(0)\cap S^3_R$.} disk centered at $0$.
Then, as in the proof of lemma \ref{lem:chi}, by additivity of
the Euler characteristic we have $\chi(\F_0) = \chi(T_0) = 1$.
Since the generic fiber is connected then $T_0$ is connected and
$H_1(T_0)= \zero$.
The morphism $j_0 : H_1(\F_0)\longrightarrow H_1(T_0)$ induced by inclusion,
is an isomorphism if and only if $0$ is a regular value at infinity
(see \cite{ACD} for the case where $\F_0$ is connected, and \cite{B}
for the general case). In our situation $j_0$ is an isomorphism
since $H_1(\F_0)= H_1(T_0)=\zero$ hence $\Binf = \varnothing$.

Conversely, let suppose now $\Baff \subset \zero$ and $\Binf =\varnothing$.
As $\Binf = \varnothing$ then $\F_0=f^{-1}(0)$ has the homotopy type of
$f^{-1}(\Delta) \cap B^4_R$. But, always because  there is no critical value at infinity, 
$f^{-1}(\Delta) \setminus B^4_R$ is just a product
$\left( f^{-1}(\Delta)\cap S^3_R \right) \times ]0,+\infty[$ (where
$f^{-1}(\Delta)\cap S^3_R$ is a tubular neighborhood in $S^3_R$ of
the link $f^{-1}(0)\cap S^3_R$). Then 
$f^{-1}(\Delta) \cap B^4_R$ has the same homotopy type as $f^{-1}(\Delta)$.
Now the polynomial $f$ is primitive since $f$ is reduced and $\B \subset \{0\}$,
hence   $f^{-1}(\Delta)$ is connected. So $\F_0$ is a connected set.
As the Euler characteristic of the connected set $\F_0$ is $+1$, it implies
that all the irreducible components of $\F_0$ are disks (possibly singular),
crossing together, without cycle (lemma \ref{lem:nocycle}). 
As a conclusion $\F_0$ is a simply connected set.
\end{proof}

\begin{remark*}
As a corollary if $\B = \zero$ with $\Binf = \zero$ then
$\F_0$ is not connected: by contraposition if $\F_0$ is a connected set then,
as $\chi(\F_0) = +1$, all irreducible components are disks (possibly singular).
Since there is no cycle (lemma \ref{lem:nocycle}) then
$\F_0$ is simply connected, thus $\Binf = \varnothing$.
\end{remark*}

\bigskip

Za\u{\i}denberg-Lin theorem admits the following generalization  
when $f$ is not irreducible.
\begin{theorem}
\label{th:ZLg}
Let $f$ be a reduced polynomial with $\Baff = \zero$ and 
$\Binf = \varnothing$ then 
$$f\sim x g(y) \qquad \text{ or } \qquad 
f \sim x^\epsilon y^{\epsilon'}\prod_{i=1}^n\big({x^p-\alpha_i y^q}\big),$$
with a non-constant polynomial $g \in \Cc[y]$, 
$\epsilon,\epsilon' \in \{0,1\}$,  $p, q$ relatively prime
numbers, $n \geqslant 1$, and $\{\alpha_i\}_{i=1,\ldots,n}$ a family of 
distinct non-zero complex numbers.
\end{theorem}

This is stated in \cite{ZL} and proved in \cite{Z} 
(see the remark after theorem $3$ of \cite{Z}).
Our proof uses topological ideas from \cite{BF}, 
\cite{NR} and \cite{R}, particularly it uses Za\u{\i}denberg-Lin 
theorem for an irreducible component of the polynomial $f$.
We need the following lemma which is a stronger version of Abhyankar-Moh
theorem.
\begin{lemma}
\label{lem:twodisks}
Let $\C_1$ and $\C_2$ be two smooth disks with equations
$(f_1=0)$ and $(f_2=0)$. 
\begin{itemize}
  \item If  $\C_1 \cap \C_2 = \varnothing$ then 
$f_1 f_2 \sim x(x+\alpha)$, $(\alpha \in \Cc^*)$.
  \item If $\C_1$ and $\C_2$ have one transversal intersection then
$f_1 f_2 \sim xy$.
\end{itemize}
\end{lemma}

\begin{proof}
By  Abhyankar-Moh theorem an equation for  $\C_1$
is $(x=0)$. As in \cite{NR} a parameterization of  $\C_2$
is $(P(t),Q(t))$, with $P,Q \in \Cc[t]$. If $\C_2$ does 
not intersect  $\C_1$ then $P(t)$ is a non-zero constant.
If  $\C_2$ intersects  $\C_1$ transversally then,
as in the proof of  Abhyankar-Moh theorem by W.~Neumann and
L.~Rudolph in \cite{NR},  polynomial automorphisms of type
$(x,y)\mapsto (x,y+\lambda x^\mu)$ enable us to choose
$(y=0)$ as an equation of $\C_2$.
\end{proof}

\begin{proof}[Proof of the theorem]
If there is no essential singularity, then singularities are ordinary quadratic
singularities. As $\F_0$ is  simply connected then it is the union of
smooth disks $\C_1,\ldots,\C_r$.
Let us suppose that $\C_1$ and $\C_2$ intersect
transversally. Then, by the lemma above, an equation of
$\C_1 \cup \C_2$ is $(xy=0)$, moreover another disk $\C_3$  can not intersect 
 $\C_1$ and $\C_2$ otherwise there is a cycle in
$G_0$ or an essential singularity. Then $\C_3$ has
equation, for instance,  $(y+\beta_3=0)$. The other
disks $\C_i$, $i\geqslant 4$ are parallel to $\C_3$ otherwise
there are cycles in $G_0$, thus $\C_i$ has equation $(y+\beta_i=0)$.
Then $f$ is algebraically equivalent to $xy\prod_i(y+\beta_i)$.

\medskip

Let us suppose that there is an essential singularity, then by lemma
\ref{lem:tubseif} there is only one essential singularity. All the other 
singularities are ordinary quadratic singularities.
Moreover, as $\Binf = \varnothing$ and as the tube $f^{-1}(S^1_r(0))$ is a Seifert
manifold (lemma \ref{lem:tubseif}), then the link at infinity 
$f^{-1}(0)\cap S^3_R$  ($S^3_R$ is a $3$-dimensional sphere with radius $R\gg 1$) 
is a Seifert link (that is to say $S^3_R\setminus f^{-1}(0)$ admits a Seifert fibration,
and the components of the link are fibers for this fibration).
By \cite[th.~2.7]{NR}, as $\Baff = \zero $ and $\Binf = \varnothing$, 
this link at infinity is the connected sum of the local
links of the singularities of $f^{-1}(0)$, that is to say
the link at infinity is the connected sum of the local link
of the essential singularity with Hopf links. 
But a Seifert link can not have such a structure, then there is only
one singularity. So the local link  and the link at infinity are
isotopic and are a sublink of
$$\OO_1 \cup \OO_2 \cup\OO(p,q)\cup\OO(p,q)\cup\ldots$$

Let explain the notations: in the sphere $S^3_r$ of $\Cc^2$, $\OO_1$, $\OO_2$
are unknots such that $\OO_1\cup\OO_2$ is the Hopf link;
 $\OO(p,q)$ denotes a torus knot of type $(p,q)$ 
($p$ and $q$ are relatively prime non-zero natural numbers)
such that $\OO_1 \cup \OO_2 \cup\OO(p,q)\cup\OO(p,q)\cup\ldots$ is isotopic to the link
$$\big((x=0)\cup (y=0) \cup (x^p-y^q=0) \cup (x^p-2y^q=0) \cup \ldots \big) \cap S^3_r.$$

\medskip

We now prove that $f$ can be written, up to algebraic equivalence,
as required. We discuss according to the number of smooth disks
in $\F_0$.
\paragraph{First case.} 
We assume that, in $\F_0$, there are two smooth 
disks with transversal intersection at the essential singularity.
By lemma \ref{lem:twodisks}, up to algebraic equivalence, an equation
of $f$ is $xy g_1(x,y)\ldots g_n(x,y)$.
We have to prove that $g_i(x,y) = x^p-\alpha_iy^q$.
Let us consider the polynomial $xyg_i(x,y)$ and
let $(P(t),Q(t))$ be  a polynomial injective parameterization  of the curve
$(g_i(x,y)=0)$.
The local link for the locally irreducible singularity is a link of type $\OO(p,q)$ 
so this parameterization can be written:
$$\begin{cases}
P(t) =  a_q t^q + \cdots + a_{N} t^{N}& \\
Q(t) =  b_p t^p + \cdots + b_{M} t^{M}& \\
\end{cases}$$
with $N \geq q$ and $M \geq p$.
As $(0,0)$ is the only point of intersection between $(xy=0)$ and $(g_i(x,y)=0)$
then  $P(t)=0$ implies $Q(t)=0$ and then $t=0$. So $P$ is  monomial: 
$P(t)=a_pt^p$. For similar reasons $Q(t) = b_qt^q$, and then $g_i(x,y)= x^p-\alpha_i y^q$.

\paragraph{Second case.} 
If $\F_0$ has only one smooth disk then for some coordinates an equation of
$f$ is $x g_1(x,y)\ldots g_n(x,y)$. As before we denote by $(P(t),Q(t))$  a parameterization
of $(g_i(x,y)=0)$; we obtain again $P(t) = a_qt^q$ $(q\ge 2)$ but
with $Q(t) = b_p t^p + \cdots + b_{M} t^{M}$. We can conclude as in 
\cite[p.~434]{NR}: the parameterization $(a_qt^q,Q(t))$ is injective so
$Q(t)-Q(\zeta t)$ has only one root $t=0$ for all $q$-th root $\zeta$ of unity. 
Hence $Q(t)$ is of the form $Q(t) = b_pt^p +F(t^q)$ and the polynomial change of
coordinates $(x,y)\mapsto (x,y - F(x))$, that preserves the axis $(x=0)$, gives
a parameterization $(a_qt^q,b_pt^p)$ in these new coordinates. So
$g_i(x,y) = x^p-\alpha_i y^q$ in these coordinates.

\paragraph{Third case.} 
If $\F_0$ has no smooth disk we can assume
that for the decomposition $f=g_1\ldots g_n$ we have $g_1(x,y) = x^p-y^q$
by Za\u{\i}denberg-Lin theorem for the irreducible component $(g_1=0)$.
Let $(g_i=0)$ be another component with parameterization 
$(P(t),Q(t)) = (a_q t^q + \cdots + a_{N} t^{N}, b_p t^p + \cdots + b_{M} t^{M})$.
The link at infinity for $(g_i=0)$ is an iterated torus knot of type 
$\OO(m,n;m_2,n_2;\ldots;m_k,n_k)$  (see \cite{R}), with
$m = M /\gcd(M,N)$ and $n=N/\gcd(M,N)$.
But the link at infinity for $(g_i=0)$ is isotopic to local link of
the affine singularity of $(g_i=0)$ and then is of type $\OO(p,q)$. 
As in  \cite{NR} either $\OO(m,n)=\OO(p,q)$ and then 
$\gcd(M,N)=1$ so $M=p$, $N=q$ and the result is proved;
or $\OO(m,n)$ is the unknot and then  $M$ divides $N$
or $N$ divides $M$. It implies that $qM\not= pN$.
As the components $(g_1=0)$ and  $(g_i=0)$ have only one intersection
(at $(0,0)$) the one variable polynomial $g_1(P(t),Q(t))$ is equal to $t^\ell$.
For example if we assume that $qM>pN$ then $\ell = qM$. 
But the valuation of $g_1(P(t),Q(t))$ 
is the intersection multiplicity of  $g_1$ and  $g_i$ at $(0,0)$ 
and it is equal to  $pq$.
Thereby $\ell = pq$  and $M=p$ and as before $N=q$.
\end{proof}

\section{Case $\Baff = \varnothing$ and $\Binf= \zero$}
\label{sec:first}

Let denote $f= f_1\times \cdots\times f_r$ the decomposition
of $f$ into irreducible factors, let $\C_i = f_i^{-1}(0)$ be
the plane algebraic curve associated to $f_i$.
We firstly obtain an ``abstract'' classification: we describe
the $\C_i$'s as punctured spheres.
\begin{proposition}
\label{prop:first}In the case $\Baff = \varnothing$ and $\Binf= \zero$, we can reorder 
the $( \C_i )_i$ so that 
\begin{itemize}
  \item either $\C_1$ is a disk and for $i = 2,\ldots,r$, $\C_i$ is an annulus;
  \item or $\C_1,\ldots,\C_{r-1}$ are disks and $\C_r$ is a $r$-punctured sphere.
\end{itemize}
\end{proposition}
\noindent
This proposition has been obtained independently in \cite{GP}.

\begin{proof}
Notice that, since $\Baff = \varnothing$ the components $\C_i$ ($i=1,\ldots,r$)
are disjoint, then
$$ \chi\big(\C_1\big)+ \cdots + \chi\big(\C_r\big) = \chi\big(\F_0\big) = 1,$$
and one of the component has positive Euler characteristic.
But as $\chi(\C_i) \leqslant 1$ for all $i$, we can suppose that the components of Euler 
characteristic $+1$ are $\C_1,\ldots,\C_j$ ($j\geqslant1$).

We firstly assume that $j= 1$; all the other components verify 
 $\chi(\C_i) \leqslant 0$ for $i \geqslant 2$, this implies that  
$\chi(\C_i) = 0$ ($i= 2,\ldots, r$).
As a conclusion the component $\C_1$ is a disk and the others are annuli.

Secondly we suppose that  $j \geqslant 2$. 
Because of the Abhyankar-Moh theorem (lemma \ref{lem:twodisks}) we can assume
that these disks $\C_1,\ldots,\C_j$ are parallel lines 
with equation $(x=\alpha_1),\ldots, (x=\alpha_j)$. 
All the other components $\C_i$ ($i>j$) have at least $j+1$ branches at infinity 
because of the non-intersection with the lines: such a component $\C_i$
has $j$ branches at infinity whose tangents at infinity are the $j$ parallel lines
that intersect the line at infinity at one point; if there is one other point
at infinity for $\C_i$ then there exists one other branch, if there is no other point
at infinity then the line at infinity is tangent to $\C_i$, that gives one more branch.
Particularly  we have $\chi(\C_{i}) \leqslant 2 - (j+1)$ for $i>j$; then
\begin{align*}
1 = \chi\big(\F_0\big) 
  &=  \chi\big(\C_1\big)+  \cdots + \chi\big(\C_j\big)
 + \chi\big(\C_{j+1}\big) + \cdots + \chi\big(\C_{r}\big) \\
  &\leqslant    \qquad  \qquad  j \qquad\quad\  \     +(1-j) +\cdots+(1-j) \leqslant 1. 
\end{align*}
Thus this inequality is an equality; it implies that $j+1=r$ and
 $\chi(\C_{j+1}) = 2 - (j+1)=2-r$, particularly there are exactly $r$
 branches at infinity.
All the components are disks, except the last one which is a $r$-punctured sphere.
This completes the proof.
\end{proof}

\bigskip

We now need a non-compact version of the Riemann-Hurwitz formula.
 Let $\bar{\C}$ and $\bar{\C}'$ be compact
Riemann surfaces and let $\pi : \bar{\C} \longrightarrow \bar{\C}'$
be a surjective holomorphic map of degree $n$. For $\mathcal{S}$ a finite
set of points in $\bar{\C}'$ we denote 
$\C' = \bar{\C}'\setminus \mathcal{S}$
and $\C = \bar{\C} \setminus \pi^{-1}(\mathcal{S})$. For any point
$s\in \C'$, $\nu(s)$ is the multiplicity of $\pi$ at $s\in \C'$, we
have $\sum_{t\in\pi^{-1}(s)}(\nu(t)-1)=n-\#\pi^{-1}(s)$.
\begin{theorem}[Riemann-Hurwitz formula]
$$\chi(\C) = n.\chi(\C') - \sum_{s\in\C}(\nu(s)-1).$$
\end{theorem}

\begin{proof}
 The proof is similar to the standard proof, see for example
\cite{RH}.
By abuse, we also denote $\C'=\bar{\C}'\setminus \mathcal{N}(\mathcal{S})$
and $\C =\bar{\C}\setminus \pi^{-1}(\mathcal{N}(\mathcal{S}))$ where 
$\mathcal{N}(\mathcal{S})$ is the union of small open disks 
around the points of $\mathcal{S}$. Let $(V',E',F')$ be a triangulation of
$\C'$ with ramification points contained in $\pi^{-1}(V')$, we denote 
$v' = \# V'$, $e'=\# E'$,...
There exists a triangulation  $(V,E,F)$ of $\C$ above $(V',E',F')$ such that
$e=ne'$, $f=nf'$ and $v = nv'-\sum_{t\in V'}(n -\# \pi^{-1}(t))$.
Then $\chi(\C)=f-e+v=n(f'-e'+v')-\sum_{t\in V'}(n- \# \pi^{-1}(t))
= n\chi(\C')-\sum_{s\in\C}(\nu(s)-1)$.
\end{proof}

We will use this formula for a component $\C=\C_i$ that is not a disk
with the natural compactification $\bar{\C}$ of $\C$:
 $\bar{\C}= \C\cup\C_\infty$.
We define $\bar{\C}' = \Cc P^1$ and if
the disks $\C_1,\ldots,\C_j$ have equation $(x=\alpha_1),\ldots,(x=\alpha_j)$
we set $\mathcal{S}= \{\infty,\alpha_1,\ldots,\alpha_j\}$
and define $\C' = \Cc P^1 \setminus \mathcal{S} = 
\Cc \setminus \{\alpha_1,\ldots,\alpha_j\}$.
The projection $\pi : \C \longrightarrow \C'$ is defined by
$\pi(x,y)=x$. Then $\pi$ can be continued to a holomorphic map on $\bar{\C}$.
If we prove that $\pi^{-1}(\mathcal{S})=\C_\infty$ then we can apply 
Riemann-Hurwitz formula.

\bigskip

We can give the algebraic classification.
\begin{proposition}
\label{prop:alg}
Depending on the cases of  proposition \ref{prop:first} above,
for a reduced polynomial $f$ with $\Baff = \varnothing$ and $\Binf = \zero$ then
\begin{itemize}
  \item either 
${f \sim x \sigma^\epsilon\prod_{i=1}^n\left({x^p\sigma^q-\alpha_i}\right)}$ or
${f \sim x \sigma^\epsilon\prod_{i=1}^n\left({ x^p-\alpha_i\sigma^q}\right)}$ 
with $p$ and $q$ relatively prime, $\epsilon \in \{ 0,1\}$, $n\geqslant 1$, 
$\{ \alpha_i \}_{i=1,\ldots,n}$ a finite family of distinct
non-zero complex numbers. Moreover
 $\sigma = \sigma(x,y) = x^sy+\ell(x)$, with $s\geqslant 0$, $\ell \in \Cc[x]$ and 
$\deg \ell < s$  (if $s>0$ then $\ell(0)\not=0$, if $s=0$ then $\ell=0$). If
$\epsilon = 1$ then $s>0$ and then in the first case $p$ (or $q$) is greater
than $1$.
  \item or ${f \sim g_\red(x)\big( g(x)y+h(x) \big)}$ with
$g,h \in \Cc[x]$, $\deg g \geqslant 2$, $\deg g > \deg h$ and $h(t) \not= 0$
if $g(t)=0$.
\end{itemize}
\end{proposition}

The situation of the first case of proposition \ref{prop:first}
has been studied by S.~Kaliman, we sketch the beautiful proof
of \cite{K}. Let $g$ be an equation of the algebraic curve $\C_1\cup
\C$ where $\C_1$ is a smooth disk and $\C=\C_i$  ($i\geqslant2$) is a 
disjoint annulus.

\begin{lemma}
$$g\sim x\big(x^p-\sigma^q\big)\qquad\text{ or }\qquad g\sim x\big(x^p\sigma^q-1\big)$$
$p,q$ are relatively prime natural numbers, and 
$\sigma=\sigma(x,y)=x^sy+\ell(x)$
($s>0$ in the first polynomial)
with $\ell\in \Cc[x]$, $\deg \ell <s$ and $\ell(0) \not= 0$ if $s>0$.
\end{lemma}

\begin{proof}
By Abhyankar-Moh theorem we can assume that $(x=0)$ is the equation for
the disk $\C_1$, let $k(x,y)$ be an equation of $\C$ in these coordinates,
there exists $m>0$ such that $k(x,x^{-m}y)= x^eh(x,y)$ with
$e < 0$, $h \in \Cc[x,y]$ and $h(0,y)=y^n$, $n \geqslant 1$.

If $\C'$ denotes the curve of equation $(h=0)$ then the ``blow-up''
$(x,y) \mapsto (x,x^{-m}y)$ gives an isomorphism from $\C'\setminus \{(0,0)\}$ 
to $\C$, 
so $\C'$ is  homeomorphic to a disk and according to 
Za\u{\i}denberg-Lin theorem the polynomial $xh(x,y)$ equal
$u(u^{p'}-v^{q'})$, the new coordinates are given by $u=x$ and
$v= y+\varphi(x)$, then $h(x,y)=x^{p'}-(y+\varphi(x))^{q'}$. Returning to $k$
by $k(x,y)=x^e h(x,x^my)$, and distinguishing the cases $e+p'=0$ and $e+p'>0$ 
leads to $k(x,y)= 1 -x^p\sigma^q$ and $k(x,y)=x^p-\sigma^q$, with
$\sigma(x,y)= x^sy+\ell(x)$. By  triangular automorphisms 
$(x,y) \mapsto (x,y+\lambda x^\mu)$ we can assume that 
$\deg \ell  < s$. That ends the proof.
\end{proof}
The generalization to the case where there are several annuli, cor\-res\-ponds
to the generalization of Za\u{\i}denberg-Lin theorem (theorem \ref{th:ZLg}).

\begin{proof}[Proof of proposition \ref{prop:alg}]
We deal with the second case of proposition \ref{prop:alg},
the disks are given by an equation $\prod_{i=1}^{r-1}(x-\alpha_i)$
and the equation of $\C_r$ is 
$$\prod_{i=1}^{r-1}(x-\alpha_i)^{m_i}\big( a_m(x)y^m+\cdots+a_1(x)y\big)+h(x)$$
with $h(\alpha_i)\not= 0$, $m_i >0$.
The projection $\pi : \C_r \longrightarrow \C' = \Cc\setminus \{
\alpha_1,\ldots,\alpha_{r-1}\}$ given by $\pi(x,y)=x$, is of degree $m$ 
and verify the hypothesis
of our Riemann-Hurwitz formula since the points at infinity
of $\C_r$ correspond to $\alpha_1,\ldots,\alpha_{r-1}$.
As $\chi(\C_r)=\chi(\C')$ then $m=1$ and $a_m$ is a constant.
The equation of $\C_r$ is now $\prod_{i=1}^{r-1}(x-\alpha_i)^{m_i}y+h(x)$
and by some triangular automorphisms $(x,y) \mapsto (x,y+\lambda x^\mu)$ we can 
assume that
$\deg h < \deg \prod_{i=1}^{r-1}(x-\alpha_i)^{m_i}$.
\end{proof}


\section{Case $\Baff = \zero$ and $\Binf= \zero$}
\label{sec:second}

Notations are those of the previous paragraph.

\begin{proposition}
\label{prop:second}
For a reduced polynomial $f$ with $\Baff = \zero$
and $\Binf = \zero$
\begin{itemize}
  \item either $\C_1$ and $\C_2$ are disks, intersecting transversally, and
$\C_i$ ($i=3,\ldots,r$) are disjoint annuli;
  \item or $\C_1,\ldots,\C_j,\C_{j+1},\ldots,\C_{r-1}$ are disjoint disks and $\C_r$
is a $(j+1)$-punctured sphere. Moreover the curves $\C_{j+1},\ldots,\C_{r-1}$ intersect $\C_r$
transversally at one point.
\end{itemize}
The corresponding algebraic list is 
\begin{itemize}
  \item either ${f \sim xy \prod_i\big(x^py^q-\alpha_i\big)}$
with $p>1$ and $q$ relatively prime, $\alpha_i\in \Cc^*$.
  \item or  ${f \sim g_\red(x)k(x)\big( g(x)y+h(x) \big)}$ with
$g,h,k \in \Cc[x]$ ($g$ and $k$ non constant, $k$ reduced), $\deg h < \deg g$ and
if $g(t)=0$ then $h(t) \not= 0$.  
\end{itemize}

\end{proposition}

\begin{lemma}
\label{lem:disk}
One of the irreducible component of $\F_0=f^{-1}(0)$ is a smooth
disk.
\end{lemma}

We will make the distinction between ``smooth'' and ``smooth in $\F_0$'':
a smooth component is not necessarily smooth in $\F_0$, there may exist
singularities on this component coming from intersection with other components.

\begin{proof}[Proof of the lemma]
Let us recall that from lemma \ref{lem:tubseif} we know that
there is at most only one essential singularity, and  affine non-essential 
singularities are ordinary quadratic singularities.
The non-essential singularities at infinity correspond to a 
bamboo\footnote{Each component intersects at most two other components.}
for the divisor at infinity $\pi^{-1}(L_\infty) \cap {\pi}^{-1}(0)$ 
for the value $0$   which intersects 
the compactification of 
some smooth disks and another component (possibly singular) of $\F_0$, 
moreover the multiplicities of
$\bar{f}$ equal to $1$ on all the components of the bamboo.
A typical example is given by  Broughton polynomial $f(x,y)=x(xy+1)$, 
another example
is given in paragraph \ref{sec:topo}.

Let us notice that one of the components of $\F_0$ is a disk 
(possibly singular) because $\chi(\F_0)=+1$. 
We firstly suppose that  no affine singularity is essential.
Then affine singularities are ordinary quadratic
singularities and a disk of $\F_0$ is smooth because it 
can not intersect itself as there is no cycle
in $G_0$.

In a second time if there exists one essential affine singularity, it is unique 
and  singularities at infinity are non-essential.
As $\Binf \not= \varnothing$ such singularities do exist.
Then one of the disks associated to a non-essential singularity at infinity 
is smooth.
\end{proof}

\begin{proof}[Proof of  proposition \ref{prop:second}]
Let $\C_1$ be the disk  of lemma \ref{lem:disk}.
Let denote $\C_1,\ldots,\C_k$ the smooth disks parallel to $\C_1$.
According to Abhyankar-Moh theorem, we can suppose that equations
for these disks are $(x=\alpha_1),\ldots,(x=\alpha_k)$.
Let $\C$ be one  of the $\C_i$ ($i>k$) which
does not intersect one of the $\C_1,\ldots,\C_k$: 
such a $\C$ exists otherwise $\F_0$ is a connected set and as
$\Binf = \zero$ this is impossible by the remark below lemma \ref{lem:simp}. 
After reordering the disks $(\C_i)_{i=1,\ldots,k}$ we denote by
$\C_1,\ldots,\C_j$ ($1\leqslant j\leqslant k$) the disks
that do not intersect $\C$. Then as above (proposition \ref{prop:first})
$\chi(\C) = 1-j$, and $\C$ has exactly $j+1$ branches at infinity.
Components other than $\C, \C_1,\ldots,\C_j$ do not contribute to Euler
characteristic: that is to say the other disks have intersection 
 with one of the $\C, \C_1,\ldots,\C_j$ at one point and at exactly one point 
because $G_0$ has no cycle; components that are not
disks are annuli.

For $\C' = \Cc \setminus \{ \alpha_1,\ldots,\alpha_j\}$
we have $\chi(\C)=\chi(\C')$ and the points at infinity correspond to
$\alpha_1,\ldots,\alpha_j$. 
The Riemann-Hurwitz formula  
for the covering $\pi : \C \longrightarrow \C'$
defined by $(x,y)\mapsto x$ proves that $\pi$ is non-branched and that $\C$ is smooth.

Singularities coming from intersections with other components can only be
transversal intersections of a smooth disk and another compo\-nent: 
the key point is the topology of $\F_0$. 
First of all, to keep $\chi(\F_0) =+1$, two components with non-positive 
Euler characteristic can not intersect; in a second time a disk that intersects
the disk $\C_1$ is smooth, otherwise it contradicts the configuration for
non-essential singularities
at infinity; and finally to avoid cycles in $G_0$, only
two directions for disks (for example $(x=0)$ and $(y=0)$) can occur, hence
there are no multiple points of order greater than $2$. We have just proved
that the affine singularities were ordinary quadratic singularities.

\bigskip

We end the classification as in proposition \ref{prop:alg}.
The main difference comes from some lines
that give ordinary quadratic singularities.
If $j=1$ and $\C_1$ is not smooth in $\F_0$ then, for Euler characteristic reasons,
there can be only one more disk $\C_2$, and $\C_2$ intersects transversally $\C_1$
at one point.
With the algebraic classification of annuli, 
we see that only one kind of annuli can occur:
$$f\sim xy\prod_{i}\big(x^py^q-\alpha_i\big),$$
where $\{\alpha_i\}$ is a family of distinct non-zero complex numbers.

For similar reasons, if $j=1$ and the disk $\C_1$ is smooth in $\F_0$, then
only one annulus can occur, but disks $\C_i$ ($i=2,\ldots,r-1$) parallel
to $\C_1$ can intersect
this annulus. Then
$$f\sim x\prod_{i}\big( x - \beta_i \big) \big( x^sy+\ell(x) \big).$$
The case $j\geqslant 2$ is treated as in proposition \ref{prop:alg}
with parallel lines  added:
$$f \sim \prod_i(x-\alpha_i)\prod_i(x-\beta_i)
\Big(   \prod_i(x-\alpha_i)^{m_i} y+h(x)\Big).$$
This completes the proof.
\end{proof}

\bigskip

The tabular summarizes the algebraic list of reduced poly\-nomials with
one critical value, notations are those of theorem \ref{th:ZLg},
propositions \ref{prop:alg} and \ref{prop:second}.
\medskip
\begin{center}
\begin{tabular}{|c|>{\centering}m{11em}|>{\centering}m{11em}|} \cline{2-3}
 \multicolumn{1}{c|}\  & $\Baff = \varnothing$ & 
          $\Baff = \{ 0 \} $ \tabularnewline \hline
$\Binf = \varnothing$ & 
 $x$ & $yg(x)$  or
         ${x^\epsilon y^{\epsilon'}\prod_i({x^p-\alpha_i y^q})}$ 
                             \tabularnewline \hline
$\Binf = \{0\}$ & 
$x \sigma^\epsilon \prod_i\left({x^p\sigma^q-\alpha_i}\right)$ or 
$x\sigma^\epsilon \prod_i\left({ x^p-\alpha_i\sigma^q}\right)$ or 
$ g_\red(x)( g(x)y+h(x) )$ & 
$xy \prod_i(x^py^q-\alpha_i)$ or 
$g_\red(x)k(x)( g(x)y+h(x) ) $ \tabularnewline \hline
\end{tabular}
\end{center}


\section{Topological classification}
\label{sec:topo}

Recall that two polynomials $f$ and $g$ are \emph{topologically equivalent}
($f \approx g$) if there exists homeomorphisms $\Phi$ and $\Psi$
such that the following diagram commutes:
$$
\xymatrix{
{}\Cc^2\ar[r]^-\Phi\ar[d]_-f  & \Cc^2\ar[d]^-g \\
\Cc \ar[r]_-\Psi        & \Cc.
}
$$

Two algebraically equivalent polynomials are topologically equivalent
but the converse is false.
For example $f(x,y) = x(x^2y+1)$ and $g(x,y)=x(x^2y+x+1)$ are
topologically equivalent (they have the same colored graph, see below) 
but are not algebraically equivalent
(an algorithm to determine if two polynomials are algebraically equivalent
is given in \cite{W}).

\bigskip

For a polynomial $f$ with resolution map $\phi$, we define the 
\emph{colored graph} $G_f$.
A vertex of the dual graph of the resolution of $f$ is colored by the value of
$\phi$ on the irreducible component associated to the vertex. In the
case $\B = \{ 0 \}$ the colors are $\infty$ (that corresponds to the 
subgraph $G_\infty$),
$0$ (for $G_0$) and $\Cc P^1$ for the dicritical components. Moreover the vertices
are weighted by the auto-intersection of the component.
In our situation all the components are rational and we do not need
to add the genus for each component.

For example, here is the graph for the polynomials  $f$ and $g$ defined above.
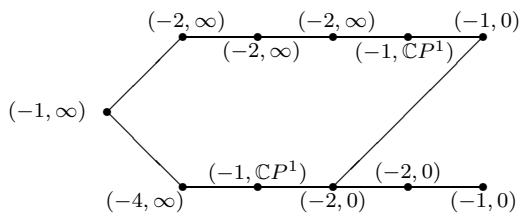
\begin{center}
\unitlength 1mm
\begin{picture}(100,30)(0,5)
\put(20,20){\circle*{1}}
\put(30,10){\circle*{1}}\put(40,10){\circle*{1}}\put(50,10){\circle*{1}}
\put(60,10){\circle*{1}}\put(70,10){\circle*{1}}
\put(30,30){\circle*{1}}\put(40,30){\circle*{1}}\put(50,30){\circle*{1}}
\put(60,30){\circle*{1}}\put(70,30){\circle*{1}}

\put(20,20){\line(1,1){10}}\put(20,20){\line(1,-1){10}}
\put(30,10){\line(1,0){10}}\put(40,10){\line(1,0){10}}\put(50,10){\line(1,0){10}}
\put(60,10){\line(1,0){10}}
\put(30,30){\line(1,0){10}}\put(40,30){\line(1,0){10}}\put(50,30){\line(1,0){10}}
\put(60,30){\line(1,0){10}}\put(70,30){\line(-1,-1){20}}

{\scriptsize
\put(12,20){\makebox(0,0){$(-1,\infty)$}}
\put(25,8){\makebox(0,0){$(-4,\infty)$}}
\put(40,12){\makebox(0,0){$(-1,\Cc P^1)$}}
\put(50,8){\makebox(0,0){$(-2,0)$}}
\put(60,12){\makebox(0,0){$(-2,0)$}}
\put(70,8){\makebox(0,0){$(-1,0)$}}
\put(30,32){\makebox(0,0){ $(-2,\infty)$}}
\put(40,28){\makebox(0,0){ $(-2,\infty)$}}
\put(50,32){\makebox(0,0){ $(-2,\infty)$}}
\put(59,28){\makebox(0,0){ $(-1,\Cc P^1)$}}
\put(70,32){\makebox(0,0){ $(-1,0)$}}
}
\end{picture}
\end{center}

Two colored graphs are \emph{equivalent} if after a sequence  of 
absorptions and blowing-ups
(see the picture below) they are isomorphic (with respect to the colors 
and weights). 
We do not authorize dicritical components to disappear in this sequence,
that is to say the color $c$ is in $\{0, \infty\}$.
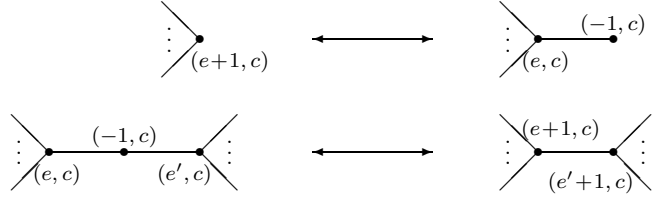
\begin{center}
\unitlength 1mm
\begin{picture}(100,25)(0,10)
{\scriptsize
\put(30,30){\circle*{1}}
\put(30,30){\line(-1,1){5}}\put(30,30){\line(-1,-1){5}}
\put(26,31){\makebox(0,0){$\vdots$}}
\put(34,27){\makebox(0,0){$(e\!+\!1,c)$}}

\put(75,30){\circle*{1}} \put(85,30){\circle*{1}} 
\put(76,27){\makebox(0,0){$(e,c)$}}\put(85,32){\makebox(0,0){$(-1,c)$}}
\put(75,30){\line(1,0){10}}
\put(75,30){\line(-1,1){5}}\put(75,30){\line(-1,-1){5}}
\put(71,31){\makebox(0,0){$\vdots$}}
\put(53,30){\vector(1,0){8}}\put(53,30){\vector(-1,0){8}}


\put(10,15){\circle*{1}}\put(20,15){\circle*{1}} \put(30,15){\circle*{1}} 
\put(11,12){\makebox(0,0){$(e,c)$}}\put(28,12){\makebox(0,0){$(e',c)$}}
\put(20,17){\makebox(0,0){$(-1,c)$}}
\put(10,15){\line(1,0){10}}\put(20,15){\line(1,0){10}}
\put(10,15){\line(-1,1){5}}\put(10,15){\line(-1,-1){5}}
\put(6,16){\makebox(0,0){$\vdots$}}
\put(30,15){\line(1,1){5}}\put(30,15){\line(1,-1){5}}
\put(34,16){\makebox(0,0){$\vdots$}}

\put(75,15){\circle*{1}} \put(85,15){\circle*{1}} 
\put(78,18){\makebox(0,0){$(e\!+\!1,c)$}}\put(82,11){\makebox(0,0){$(e'\!+\!1,c)$}}
\put(75,15){\line(1,0){10}}
\put(75,15){\line(-1,1){5}}\put(75,15){\line(-1,-1){5}}
\put(71,16){\makebox(0,0){$\vdots$}}
\put(85,15){\line(1,1){5}}\put(85,15){\line(1,-1){5}}
\put(89,16){\makebox(0,0){$\vdots$}}
\put(53,15){\vector(1,0){8}}\put(53,15){\vector(-1,0){8}}
}
\end{picture}
\end{center}

\begin{proposition}
Let two reduced polynomials have only one critical value.
If they have equivalent colored graphs then they are topologically equivalent
\end{proposition}

This proposition can not be generalized to the case of several critical values, 
a counter-example is given in \cite{AB}.
The converse is true: the main ideas for proving this
are in the proof of the next proposition or refer to \cite{F}.

\begin{proof}
Let $f$ and $g$ be polynomials with just one critical value $0$ and with  
equivalent colored graphs. Let $(\pi_f,\bar{f})$, $(\pi_g,\bar{g})$
come from the resolution of $f$ and $g$. One can suppose, after some blowing-ups and 
absorptions, that their graphs are equal. We set 
$D_{f,0}= \bar{f}^{-1}(0)$ and $D_{g,0}= \bar{g}^{-1}(0)$.

By standard arguments (\cite{AC}, \cite{Du}, \cite{F}), a small
neighborhood of $D_{f,0}$ is homeomorphic to a small
neighborhood of $D_{g,0}$. As all the components of $D_{f,0}$
and $D_{g,0}$ are rational the monodromies for $\bar{f}$ 
and $\bar{g}$ induced by a small circle 
around the value $0$ act equivalently: that is to say the following
diagram commutes:
$$
\xymatrix{
{}\bar{f}^{-1}(\Delta)\ar[r]^-{\bar{\Phi}_0}\ar[d]_-f  
               & \bar{g}^{-1}(\Delta')\ar[d]^-g \\
\Delta  \ar[r]_-{\Psi_0}        & \Delta'
}
$$
where $\bar{\Phi}_0$ and $\Psi_0$ are homeomorphisms and
$\Delta$ and $\Delta'$ topological
closed disks of $\Cc$ with $0 \in \Int \Delta \cap \Int \Delta'$.

Let $D^\infty_{f,0} = D_{f,0} \cap \pi_f^{-1}(L_\infty)$ be
the part of $D_{f,0}$ that corresponds to the irregularity at infinity of 
the value $0$ ($D^\infty_{g,0}$ is set in the same way).
Then $\bar{\Phi}_0$  defines an homeomorphism between 
$D^\infty_{f,0}$ and $D^\infty_{g,0}$. Then the homeomorphism
$\bar{\Phi}_0$ from $\bar{f}^{-1}(\Delta) \setminus 
D^\infty_{f,0}$ to $\bar{g}^{-1}(\Delta') \setminus 
D^\infty_{g,0}$ can be restricted to an homeomorphism
$\Phi_0$ that respects the fibration because $f\circ \pi_f =
\bar{f}$ on the set 
$\bar{f}^{-1}(\Delta) \setminus D^\infty_{f,0}$.
We have proved that $f$ and $g$ are topologically
equivalent in a neighborhood of the zero fiber:
$$
\xymatrix{
{}f^{-1}(\Delta)\ar[r]^-{\Phi_0}\ar[d]_-f  & g^{-1}(\Delta')\ar[d]^-g \\
\Delta  \ar[r]_-{\Psi_0}        & \Delta'.
}
$$

\bigskip

We now explain how to continue theses homeomorphisms.
As the only critical value for $f$ is in $\Delta$ the 
fibration $f : f^{-1}(\overline{\Cc \setminus  \Delta})
\longrightarrow \overline{\Cc \setminus  \Delta}$ is isomorphic
to the fibration $f\times\id : f^{-1}(\partial \Delta)\times \Rr_+ \longrightarrow 
\partial \Delta \times \Rr_+$. It provides homeomorphisms  $\phi_f$
and  $\psi_f$ (see diagrams). For $g$ we obtain homeomorphisms $\phi_g$
and  $\psi_g$. But the fibration $f\times\id$ above $\partial \Delta \times \Rr_+$
is isomorphic to the fibration $g\times\id$ above $\partial \Delta' \times \Rr_+$,
the corresponding homeomorphisms are $\Phi_1 = \Phi_0\times \id$ and 
 $\Psi_1 = \Psi_0\times \id$.
{\small 
$$
\xymatrix{
{}
f^{-1}(\overline{\Cc \setminus  \Delta})\ar[r]^-{\phi_f}\ar[d]_-f &
f^{-1}(\partial \Delta)\times \Rr_+ \ar[r]^-{\Phi_1}\ar[d]_-{f\times \id} &
g^{-1}(\partial \Delta')\times \Rr_+ \ar[d]_-{g\times \id} & 
g^{-1}(\overline{\Cc \setminus  \Delta'})\ar[l]_-{\phi_g}\ar[d]^-g \\ 
\overline{\Cc \setminus  \Delta} \ar[r]_-{\psi_f}&
\partial \Delta\times \Rr_+ \ar[r]_-{\Psi_1} &
\partial \Delta'\times \Rr_+  &
\overline{\Cc \setminus  \Delta'}. \ar[l]^-{\psi_g}\\
}
$$
}
Then $\Phi_0$ can be continued by ${\phi}_{g}^{-1} \circ
{\Phi}_{1}\circ{\phi}_{f}$ and
$\Psi_0$ by ${\psi}_{g}^{-1} \circ
{\Psi}_{1}\circ {\psi}_{f}$.
\end{proof}

\bigskip

\begin{theorem}
\label{prop:top}
A reduced polynomial with at most one critical value 
is topologically equivalent to one, and only one, of the following polynomials
(notations are those of the introduction):
\medskip
\begin{center}
\begin{tabular}{|c|>{\centering}m{11em}|>{\centering}m{10.5em}|} \cline{2-3}
 \multicolumn{1}{c|}\  & $\Baff = \varnothing$ & 
          $\Baff = \{ 0 \} $ \tabularnewline \hline
$\Binf = \varnothing$ & 
 $x$ & 
 $yg_\red(x)$  or
 $x\prod_{i=1}^{n}({x^p-i y})$ 
     or    ${x^\epsilon y^{\epsilon'}\prod_{i=1}^{n}({x^p-i y^q})}$ 
                             \tabularnewline \hline
$\Binf = \{0\}$ & 
$x \prod_{i=1}^{n}\left({x^py^q-i}\right)$ or 
$x \sigma\prod_{i=1}^{n}\left({x^p\sigma^q-i}\right)$ or 
$x\sigma^\epsilon \prod_{i=1}^{n}\left({ x^p-i\sigma^q}\right)$ or 
$ g_\red(x)( g(x)y+1)$ & 
$xy \prod_{i=1}^{n}(x^py^q-i)$ or 
$g_\red(x)k(x)( g(x)y+1 ) $ \tabularnewline \hline
\end{tabular}
\end{center}
\medskip
\noindent
\end{theorem}
\begin{proof}
We firstly  have to prove that the list of polynomials up to algebraic 
equivalence can be 
reduced, up to topological equivalence, to the list above. Finally we shall prove
that two distinct polynomials of this list are not topologically equivalent.

For the cases with $\Binf = \varnothing$, replacing $\alpha_i$ by $i$ does not
change the polynomial, up to topological equivalence. Moreover the list, for
theses cases, is not redundant.

Let study what happens for the case $\Binf = \zero$.
Let $f$ be one of the polynomials coming from the algebraic list, and 
let  $f'$ be the corres\-ponding polynomial with the constant $1$
instead of the polynomial $\ell(x)$ or $h(x)$ and with
$i$ instead of $\alpha_i$.
We may find $f \approx f'$ by proving that the graphs $G_f$  and $G_{f'}$
are equivalent.
As $f$ and $f'$  have the same behavior at finite distance,
we just have to study what happens at infinity.

Let $F(x,y,z)$  be the homogeneous polynomial
associated to $f$, $P_1=(1\!:\!0\!:\!0)$ and  $P_2=(0\!:\!1\!:\!0)$ are 
the two points 
at infinity of $f$; we denote $f_1(y,z) = F(1,y,z)$,
$f_2(x,z) = F(x,1,z)$ the local equations of $F$ at the points $P_1$,
$P_2$. 
To calculate the part of $G_f$ at infinity, we have two ---equivalent---
choices.

Firstly we can calculate the irregular link
at infinity $f^{-1}(0)\cap S^3_R$ (\emph{resp.}\ ${f'}^{-1}(0)\cap S^3_R$),
it is a sufficient condition since the (single) irregular link determines 
the regular links at infinity
$f^{-1}(s)\cap S^3_R$  ($s\not=0$),  see \cite{NT}.

Secondly, we can calculate the  
 Puiseux expansions of the branches of $f_1$ (and $f_2$)
and  the intersection multiplicities between the branches of
$f_1$ (and between the branches of $f_2$) by taking into account the
line at infinity with local equation $(z=0)$.
It is a sufficient condition since if we know
the topology of $zf_i$ then one can recover the topology
of the family $(f_i-tz^d)_{t\in \Cc P^1}$ (see \cite{LW2}) as $t=0$ and
$t=\infty$ are the only critical values for this family.

We will use the second method:
 $f \approx f'$ if and only if $f_1$ and $f'_1$ (and $f_2$ and $f'_2$) 
have equivalent  Puiseux expansions and the same intersection multiplicities.

\bigskip

We will detail the calculus for $f(x,y) = x \sigma \prod_{i=1}^{n}\left({x^p\sigma^q-\alpha_i}\right)$
with $\sigma(x,y) =x^sy+\ell(x) = x^sy+a_{s-1}x^{s-1}+\cdots+a_0$,
$a_0 \not= 0$ and $n > 1$, the calculus are similar
for the other polynomials. Then
$f'(x,y) = x {\sigma'} \prod_i\left({x^p{\sigma'}^q-i}\right)$
and $\sigma'(x,y) = x^sy+1$.
The local equation of $F$ at $P_1$ is
\begin{align*}
f_1(y,z) =&   (y+a_{s-1}z^{2}+\cdots+a_0z^{s+1}) \times \\
          &\quad  \prod_i\left({(y+a_{s-1}z^{2}+\cdots+a_0z^{s+1})^q
 - \alpha_iz^{p+q(s+1)}}\right).
\end{align*}
A similar formula holds for $f_1'$.
The branches of $f_1$ and $f'_1$ are smooth and intersect
the line at infinity $(z=0)$ transversally. Moreover the intersection
multiplicities for the branches of $f_1$ are independent of the
coefficients $a_{s-1},\ldots,a_1$, of $a_0 \not=0$, and of the $\alpha_i \not= 0$:
let $\ell_1(y,z) = y+a_{s-1}z^{2}+\cdots+a_0z^{s+1}$
then $\m(\ell_1(y,z),\ell_1^q(y,z)-\alpha_iz^{p+q(s+1)})=
p+q(s+1)$; for $i\not= j$, $\m(\ell_1^q(y,z)-\alpha_jz^{p+q(s+1)},
\ell_1^q(y,z)-\alpha_iz^{p+q(s+1)})=q(p+q(s+1))$
(see how to calculate intersection multiplicities below),
so $f_1$ and $f'_1$ have equivalent Puiseux expansions and 
the same intersection multiplicities.

The following lemma allows us to calculate intersection multiplicities;
the first point is well-known (see \cite{BK} or \cite{D}), 
the second point is a consequence of the first.
\begin{lemma}Let $f,g,f_1,f_2$ be irreducible plane curve germs at $0$.
\begin{itemize}
\item
 Let $K_f,K_g$
be the local links of $f$ and $g$. Then the intersection multiplicity verify
\begin{align*}
\m(f,g) &=\dim_\Cc \Cc \lbrace x,y \rbrace /(f(x,y),g(x,y))\\ 
        &= \mathrm{lk} (K_f,K_g)
        =\mathrm{val}_t (f\circ p(t))
\end{align*}
with $\mathrm{lk}$ is the linking number, $\mathrm{val}$ is the valuation
and $p(t) = (t^n, \varphi(t))$ is a Puiseux parameterization for the
curve $(g=0)$ (which is supposed not to contain $(y=0)$).
\item Let $t,t'$ be complex numbers with $t\not= t'$ and $t'\not=0$.
Then
\begin{align*}
&\m(f_1f_2,g) = \m(f_1,g)+\m(f_2,g), \\
&\m(f+tg,f+t'g)=\m(f,f+t'g).
\end{align*}
\end{itemize}
\end{lemma}

\noindent
For the second point at infinity $P_2$
the local equation of $F$ is 
\begin{align*} 
f_2(x,z) &=  x (x^s+a_{s-1}x^{s-1}z^{2}+\cdots+a_0z^{s+1}) \times \\
 &\quad \prod_i\left({x^p(x^s+a_{s-1}x^{s-1}z^{2}+\cdots+a_0z^{s+1})^q
 - \alpha_iz^{p+q(s+1)}}\right).
\end{align*}
All the branches intersect transversally the line at infinity, and 
the topology of each branch is given by one of the Puiseux expansions
 $x=0$, $x= z^{\frac{s+1}{s}}$ and 
$x= z^{\frac{p+q(s+1)}{p+qs}}$
 and is independent of $a_{s-1},\ldots,a_1$, of $a_0 \not=0$ and
of the $\alpha_i \not=0$.
Moreover  intersection multiplicities are also independent of the
coefficients: let $\ell_2(x,z)=x^s+a_{s-1}x^{s-1}z^{2}+\cdots+a_0z^{s+1}$
then $\m(x,\ell_2(x,z)) = s+1$, $\m(x,x^p\ell_2(x,z)^q-\alpha_jz^{p+q(s+1)})
= p+q(s+1)$, $\m(\ell_2(x,z),x^p\ell_2(x,z)^q-\alpha_jz^{p+q(s+1)})
= s(p+q(s+1))$, and for $i \not= j$, 
$\m(x^p\ell_2(x,z)^q-\alpha_iz^{p+q(s+1)},
x^p\ell_2(x,z)^q-\alpha_jz^{p+q(s+1)})= (p+qs)(p+q(s+1))$.

As a conclusion $f_1$, $f'_1$ and $f_2$, $f'_2$ have
the same branches and the branches have the same tangency,
so $f$ and $f'$ are topologically equivalent.

\bigskip

Finally, we shall prove that the list is non-redundant.
As before, we detail the calculus for the polynomial
$f(x,y)= x \sigma \prod_{i=1}^n(x^p\sigma^q-i)$ with
$\sigma = x^sy+1$; for the other polynomials the method is the same.
Let suppose that another polynomial, $f'$, of the topological list
verify $f \approx f'$. Then $f'$ has the same type as $f$, that is to say
that $f'(x,y) = x\sigma'\prod_{i=1}^n(x^{p'}{\sigma'}^{q'}-i)$ with
$\sigma' = x^{s'}y+1$.
As $f \approx f'$, the localizations $f_1$ and $f'_1$ (\emph{resp.}\ 
 $f_2$ and $f'_2$) at $P_1$ (\emph{resp.}\  $P_2$) have equivalent
Puiseux expansions and the same intersection multiplicities. 
We deduce from the calculus of intersection multiplicities
at $P_2$ that $s+1=s'+1$ and at $P_1$ that
$p+q(s+1) =p'+q'(s'+1) $ and $q(p+q(s+1)) =q'(p'+q'(s'+1)) $.
It implies that $s=s'$, $p=p'$, $q=q'$ and then $f=f'$.
\end{proof}

\bigskip

{\small
\emph{Acknowledgements:} I would like to thank Michel Boileau and Fran\c{c}oise Michel for their encouragements.
}
\medskip


{\noindent
Arnaud \textsc{Bodin} \\
Universit\'e Paul Sabatier Toulouse III, laboratoire \'Emile Picard, \\
118 route de Narbonne, 31062 Toulouse c\'edex 4, France. \\
\texttt{bodin@picard.ups-tlse.fr}
}

\end{document}